\documentclass[11pt]{amsart}
\usepackage{amsfonts,amssymb,amsmath}

\newtheorem{Theorem}{Theorem}
\newtheorem{proposition}[Theorem]{Proposition}

\newtheorem{thm}[Theorem]{Theorem}
\newtheorem{theorem}[Theorem]{Theorem}

\theoremstyle{definition}

\newcommand{\alg}{\mathbf} 
\def\v{\mathcal V}  

\newcommand{\brfr}{$\hspace{0 pt}$} 

\newcommand{\adma}{\mathrm{Adm({\alg A})}} 
\newcommand{\admb}{\mathrm{Adm({\alg B})}} 
\newcommand{\admx}{\mathrm{Adm({\alg X})}} 

\DeclareMathOperator{\RR}{\mathit{R}} 
\DeclareMathOperator{\SSS}{\mathit{S}} 
\DeclareMathOperator{\TT}{\mathit{T}} 
\DeclareMathOperator{\UU}{\mathit{U}} 
 
\DeclareMathOperator{\FT}{\mathit{F(T)}} 

\DeclareMathOperator{\GS}{\mathit{G(S)}}

\begin{document}

\title{Towards Commutator theory for relations. IV.}

\author{Paolo Lipparini}
\address{Dipartimento di Matematica,
 Viale della Ricerca Scientifica,
II Universit\`a di Roma (or VergaTta),
 ROME 
ITALY}

\thanks{The author has received support from MPI and GNSAGA} 


\urladdr{http://www.mat.uniroma2.it/\textasciitilde lipparin}

\begin{abstract} 
We provide more characterizations of 
varieties having a term Mal'cev modulo two functions $F$ and $G$.
We characterize varieties neutral in the sense of $F$, that is
varieties satisfying $R \subseteq F(R)$.
We present examples of global operators satisfying the homomorphism property, in particular we show that many known commutators satisfy the homomorphism property.
\end{abstract} 

\maketitle

See Parts I-III \cite{ctfr} for unexplained notation.

We now find more conditions equivalent to the existence of 
a term 
Mal'cev modulo $F$ and $G$ throughout a variety,
thus complementing Theorem 3 in Part III.

\begin{thm} \label{wdtvarIVa}
Suppose that $\v$ is a variety,
 $F$, $G$ are 
global operators on $\v$ for admissible and reflexive relations, 
$F$, $G$ are monotone and satisfy the homomorphism property.
Then the following are equivalent:

(i)
$\v$ has a term which is 
Mal'cev modulo $F_{\alg A}$ and $G_{\alg A}$
for every algebra ${\alg A}$ in $\v$.

(ii) In every algebra ${\alg A}\in \v$ and 
for all 
relations $R, S, T \in \adma$, the following holds:
\[
T(R \circ S) \subseteq F_{\alg A} (T) \circ \overline{S^-\cup R} \circ G_{\alg A} (S).
\]

(iii) In the free algebra ${\alg X}$ in $\v$ generated by $2$ elements the following holds:
\[
T(R \circ S) \subseteq F_{\alg X} (T) \circ S^-\circ R \circ G_{\alg X} (S),
\]
for all 
relations $R, S, T \in \admx$.

(iv) In every algebra ${\alg A}\in \v$ and 
for all 
relations $ S, T \in \adma$, the following holds:
\[
TS \subseteq F_{\alg A} (T) \circ S^- \circ G_{\alg A} (T).
\]

(v) In the free algebra ${\alg X}$ in $\v$ generated by $2$ elements the following holds:
\[
TS \subseteq F_{\alg X} (T) \circ S^- \circ G_{\alg X} (T),
\]
for all 
relations $S, T \in \admx$.

(vi) In the free algebra ${\alg X}$ in $\v$ generated by the $3$ elements 
$ x,y,z $ the following holds:
\[
S \circ T \subseteq F_{\alg X} (S) \circ Cg(T) \circ Cg(S) \circ G_{\alg X} (T).
\]
where $S $ is the smallest reflexive and admissible relation on $ {\alg X} $
containing $(x,y)$,
and  $ T$ is the smallest reflexive and admissible relation on $ {\alg X} $
containing $(y,z)$.

(vii) In the free algebra ${\alg X}$ in $\v$ generated by the $3$ elements 
$  x,y,z $ the following holds:
\[
S  \circ S  \subseteq F_{\alg X} ( S ) \circ \Theta  \circ G_{\alg X} ( S ).
\]
where $ S $ is the smallest reflexive and admissible relation on  $ {\alg X} $
containing $(x,y)$ and $(y,z)$, and
$ \Theta  $ is the  tolerance generated by $ S$.  
\end{thm} 

\begin{proof}
(i) $\Rightarrow $ (ii). If  $(a,c) \in T(R \circ S)$
then $a \TT c$ and there exists $b \in A$ such that 
$a \RR b \SSS c$. Thus,
$a \FT t(a,c,c) \mathrel{\overline{\SSS^- \cup \RR}} t(b,b,c) \GS c$.
Actually, (i) $\Rightarrow $ (ii) is a particular case of
Theorem 1(i) in Part III.

(ii) $\Rightarrow $ (iii) and (iv) $\Rightarrow $ (v) are trivial. 

Both (ii) $\Rightarrow $ (iv) 
and
(iii) $\Rightarrow $ (v)
are obtained by taking $R=0$.

(v) $\Rightarrow $ (i). Let $x,y$
be the generators of $ {\alg X} $. By taking, in (v), $T=S$  to be
the smallest admissible and reflexive relation of ${\alg X}$ containing
$(x,y)$ we get Condition (x) in Part III, Theorem 3, which implies (i)
by the same theorem.

(i) $\Rightarrow $ (vi) 
and
(i) $\Rightarrow $ (vii) 
are immediate from Part III, Theorem 1 (iii).

(vi) $\Rightarrow $ (i).  
Since
$ x \SSS y \TT z $, we have 
$ (x,z) \in S \circ T$,
hence 
$(x,z) \in F_{\alg X} (S) \circ Cg(T) \circ Cg(S) \circ G_{\alg X} (T)$
by assumption.

This means that there are terms $t_1(x,y,z)$, $t(x,y,z)$ and  
$t_2(x,y,z)$ such that 
 $(x,
 t_1(x,y,z)) \in F_{\alg X} (S)$,
$(t_1(x,y,z), t(x,y,z)) \in Cg(T)$,
$(t(x, y, z), 
t_2(x,
\brfr
y,z)) \in Cg(S)$ and 
$( t_2(x,y,z),z) \in G_{\alg X} (T)$.

Since $ Cg(T)= \{(u(x,y,z), v(x,y,z))| u, v \text{ terms of } 
{\v} \text{ such that }  u(x,y,y)= v(x,y,y) \} $,
we get
$t_1(x,y,y)=t(x,y,y)$,
and similarly 
$t(y,y,z)=t_2(y,y,z)$.

Now the proof goes exactly as in the final part of the 
proof of Part III, Theorem 3
(vii) $\Rightarrow $ (i).  

(vii) $\Rightarrow $ (i).
Let $ {\alg X}, S, \Theta  $ be as in (vii).
Since $x \SSS y \SSS z$, we have
$(x,z) \in S  \circ S $, hence 
$ (x,z) \in  F_{\alg X} ( S ) \circ \Theta  \circ G_{\alg X} (S )$.

This means that $\v$ has terms $t_1(x,y,z)$ and  
$t_2(x,y,z)$ such that 
 $(x,
\brfr
 t_1(x,y,z)) \in F_{\alg X} ( S)$,
$(t_1(x,y,z), t_2(x,y,z)) \in \Theta$ and 
$( t_2(x,y,z),z) \in G_{\alg X} (S)$.

Notice that $ \Theta= \{(u(x,y,z,x,y,y,z), u(x,y,z,y,x,z,y))| u \text{ a term of } 
{\alg X} \} $. Indeed, the right-hand relation is  reflexive, symmetric,
admissible, and contains $(x,y)$ and $(y,z)$; moreover, every other reflexive, 
symmetric admissible
relation containing $(x,y)$ and $(y,z)$ contains all pairs of the form
$(u(x,y,z,x,y,y,z), u(x,y,z,y,x,z,y))$.

Hence, $(t_1(x,y,z), t_2(x,y,z)) \in \Theta$ means that 
there is
 a $7$-ary term $t'$ such that   
$t_1(x,y,z)=t'(x,y,z,x,y, y,z)$ and $t'(x,y,z,y,x,z,y)= t_2(x,y,z) $.

We claim that the ternary term 
$t(x,y,z)=t'(x,y,z,x,y,z,y)$ is Mal'cev modulo $F$ and $G$ throughout $\v$.
First, notice that
$t_1(x,y,y)=t(x,y,y)$, and 
$t(y,y,z)=t_2(y,y,z)$.
Now the last part of the proof of 
Part III, Theorem 3
(vii) $\Rightarrow $ (i) shows that
$a \mathrel{F_{\alg A} (R)} t_1(a,b,b)=t(a,b,b)$
holds for every algebra $ {\alg A} \in \v$
and $a \RR b$, and a similar relation holds for
$G$, thus $t$ is Mal'cev modulo $F$ and $G$.
\end{proof}

We say that 
$F$ is \emph{regular} on $\v$ if and only if 
$F$ is a
global operator on $\v$ from admissible and reflexive relations
to congruences (that is, $F_{{\alg A}}(R)$ is always supposed to be a congruence of $ {\alg A} $), 
$F$ is monotone, satisfies the homomorphism property,
and satisfies also the following \emph{quotient property}:
``for every $ {\alg A} \in \v$ and every $R \in \adma$, it happens that
$F_{\alg B }( {R / {F_{\alg A }(R)} } )=0_ {\alg B} $,
where $ {\alg B}$ is ${\alg A / {F_{\alg A }(R)} }  $''.
Here, $ R / {F_{\alg A }(R)}$ denotes $ \pi(R)$,
where $\pi$ denotes the canonical epimorphism
$: {\alg A} \to {\alg A}/ F_{\alg A }(R)$.

In the particular case when $F=G$,
we simply say that a term $t$ is 
\emph{Mal'cev modulo}
 $F$ in order to mean that
$t$ is Mal'cev modulo $F$ and $F$.
 
\begin{theorem} \label{wdtvarIVb}
Suppose that $\v$ is a variety and
 $F$ is a 
global operator which is regular on $\v$ and satisfies
$F(R)=F(R^-)$, for every admissible relation $R$.
 Then the following are equivalent: 

(i)
$\v$ has a term which is 
Mal'cev modulo $F_{\alg A}$ for every algebra ${\alg A}$ in $\v$.

(ii) For every $ {\alg A} \in V$ and $R \in \adma$, if
$F_{\alg A}(R)=0_{ {\alg A} }$ then $R$ is a congruence of $ {\alg A} $.

(iii) For every $ {\alg A} \in V$ and $R \in \adma$, if
$F_{\alg A}(R)=0_{ {\alg A} }$ then $R$ is a tolerance of $ {\alg A} $.

(iv) For every $ {\alg A} \in V$ and $R \in \adma$, if
$F_{\alg A}(R)=0_{ {\alg A} }$ then $R$ is a transitive relation of $ {\alg A} $.
 
(v) Let $ {\alg X} $ be the free algebra in $\v$ generated by $ \{ x, y\} $ and
$S $ be the smallest reflexive and admissible relation on $ {\alg X} $
containing $(x,y)$.
If $ {\alg B} = {\alg X} /F_{{\alg X} }(S)$, then
$T=S/ F_{{\alg X} }(S)$ is (a) a congruence (equivalently, (b) a tolerance) of $ {\alg B} $.

(vi) Let $ {\alg X} $ be the free algebra in $\v$ generated by $ \{ x, y, z \} $ and
$S $ be the smallest reflexive and admissible relation on $ {\alg X} $
containing $(x,y)$ and $(y,z)$.
If $ {\alg B} = {\alg X} /F_{{\alg X} }(S)$, then
$T=S/ F_{{\alg X} }(S)$ is a transitive relation of $ {\alg B} $.
\end{theorem} 

\begin{proof} 
(i) $\Rightarrow $ (ii) is immediate from Part III, Theorem 1(xi).

(ii) $\Rightarrow $ (iii) and (ii) $\Rightarrow $ (iv) are trivial.

Let $ {\alg X}, {\alg B}, S, T $ be as in (v) (respectively, as in (vi)).
By the quotient property,  $ F_{\alg B}(T)=0_{\alg B}$, thus
(ii), (iii), (iv) imply, respectively, (v)(a), (v)(b), (vi).

If $T$
is a congruence then $T$ is a tolerance, thus (v)(a) implies (v)(b).

Suppose that (v)(b) holds.
Since $T=S/ F_{{\alg X} }(S)$ is a tolerance, 
and $\pi$ is surjective, then $\pi^{-1}(T)= \pi^{-1} (S/ F_{{\alg X} }(S))= F_{{\alg X} }(S) \circ S \circ F_{{\alg X} }(S)$ is a tolerance, hence
$S \subseteq 
F_{{\alg X} }(S) \circ S \circ F_{{\alg X} }(S)=
\big(F_{{\alg X} }(S) \circ S \circ F_{{\alg X} }(S)\big)^-=
F_{{\alg X} }(S) \circ S^- \circ F_{{\alg X} }(S)$.

Thus
$S \subseteq 
F_{{\alg X} }(S) \circ S^- \circ F_{{\alg X} }(S)$,
which is
Condition (x) in Part III, Theorem 3, which implies (i).

Suppose that (vi) holds, that is, $T$ is transitive.
Then $\pi^{-1}(T)$ is transitive, hence
$S \circ S \subseteq 
\pi^{-1}(T) \circ \pi^{-1}(T)= \pi^{-1}(T)=
F_{{\alg X} }(S) \circ S \circ F_{{\alg X} }(S)$.

This is Condition (vii) in Part III, Theorem 3, which implies (i).
\end{proof}

The next theorem provides a version of our results for what might be called
\emph{neutral} varieties, that is, varieties satisfying
$ R \subseteq F(R)$ (notice that we are not necessarily assuming $F(R)
\subseteq R$).

\begin{thm} \label{neutrvar} 
Suppose that $\v$ is a variety,
 $F$ is a
global operator on $\v$ for admissible and reflexive relations, 
and $F$ is monotone and satisfies the homomorphism property.
Then the following are equivalent:

(i)
$R \subseteq F_{\alg A} (R)$
holds
for every algebra ${\alg A}$ in $\v$ and
for every 
relation $R \in \adma$.

(ii) $R \subseteq F_{\alg X} (R)$
holds
for every
admissible and reflexive relation
$R$ in
the free algebra ${\alg X}$ in $\v$ generated by $2$ elements. 

(iii) $S \subseteq F_{\alg X} (S)$
holds
in the free algebra ${\alg X} $ in $\v$ generated by the two elements
$x$, $y $, where 
$S$ is the smallest admissible and reflexive relation 
of ${\alg X}$
containing
$(x,y)$.

(iii)$'$ 
In the free algebra ${\alg X} $ in $\v$ generated by the two elements
$x$, $y $, if 
$S$ is the smallest admissible and reflexive relation 
of ${\alg X}$
containing
$(x,y)$ then
$(x,y) \in F_{\alg X} (S)$.

In case $F$ is regular, the preceding conditions are also
equivalent to the following ones:

(iv) For every algebra ${\alg A}$ in $\v$ and
for every 
relation $R \in \adma$,
if $R \not=0$ then $ F_{\alg A}(R) \not=0$.

(v) Let ${\alg X} $  be the free algebra in $\v$ generated by the two elements
$x$, $y $, and let
$S$ be the smallest admissible and reflexive relation 
of ${\alg X}$
containing
$(x,y)$.
Then $S/F_ {\alg X}(S)=0 $ in $ {\alg B}={\alg X}/F_{\alg X}(S) $.
\end{thm} 

\begin{proof}
(i) $\Rightarrow $ (ii) $\Rightarrow $ (iii) $\Rightarrow $ (iii)$'$ 
are trivial.

(iii)$'$ $\Rightarrow $ (i). 
Suppose that
${\alg A}\in \v$, 
$a, b \in  A $, 
$R \in \adma$,
and $a \RR b$.

Since ${\alg X} $ is
the free algebra  in $\v$ generated by 
$\{x,y\} $,
there is a homomorphism
$\phi: {\alg X} \to {\alg A} $
such that 
$\phi(x)=a$
and $\phi(y)=b$.

Since, 
by (iii)$'$,
$ x \mathrel{F_{\alg X} (S)} y$,
we have that
$ \phi (x) 
\mathrel{\phi ( F_{\alg X} (S))}
 \phi ( y)$,
hence, by the homomorphism property,
$a
 \mathrel{F_{\alg A}(\phi (S))}
 b$.

Since 
$S$ is the admissible and reflexive relation generated by $(x,y)$,
and $\phi$ is a homomorphism, 
then $\phi(S)$ is the admissible and reflexive relation generated by 
$(\phi(x),\phi(y))=(a,b)$,
and, since $a \RR b$, we have that 
$ \phi(S) \subseteq R$; thus, by the monotonicity of
$F_{\alg A}$, we get
$ F_{\alg A}(\phi(S)) \subseteq F_{\alg A}(R)$ and, eventually,
$a 
\mathrel{F_{\alg A}(R)}
 b$.

(i) $\Rightarrow $ (iv) is trivial.

(iv) $\Rightarrow $ (v).
 By the quotient property,
$F_ {\alg B}(  S/F_ {\alg X}(S))=0_  {\alg B}  $, thus
$ S/F_ {\alg X}(S)=0_  {\alg B} $ by (iv).

(v) $\Rightarrow $ (iii).
$ S/F_ {\alg X}(S)=0_  {\alg B} $ in $ {\alg B} $ means
$ S \subseteq F_ {\alg X}(S)$ in $ {\alg X}$.
\end{proof}


We now show that the commutators defined in
\cite{lpumi} and in Part I
satisfy the homomorphism property. Actually, we shall deal
with $n$-ary 
global operators on $\v$, that is, operators depending on $n$-variables 
(as above, each variable is intended to be
an admissible and reflexive relation).
In this general situation, the {\em homomorphism property}
reads 
$\phi(F_{\alg B}(R_1, \dots, R_n))
 \subseteq F_{\alg A}(\phi(R_1), \dots, \phi(R_n))$,
for every 
${\alg A}, {\alg  B} \in \v $,
homomorphism $\phi: {\alg B} \to {\alg A} $,
and $R_1, \dots, R_n \in \admb$.
Moreover, $F$ is {\em  monotone}
if and only if 
$R_1 \subseteq S_1$,
\dots,
$R_n \subseteq S_n$
imply
$F(R_1, \dots , R_n) \subseteq F(S_1, \dots, S_n)$.

If $\mathbf A$ is any algebra, and $ R, S$ are 
compatible and reflexive relations,
 $M( R , S )$ is defined to be the set of all {\em matrices} of the form
\[
\begin{vmatrix}
t( \bar{a}, \bar{b})     & t( \bar{a}, \bar{b}') \cr
t( \bar{a}', \bar{b})     & t( \bar{a}', \bar{b}') 
\end{vmatrix},   
\] 
where $\bar{a}, \bar{a}' \in A^h $,
$\bar{b}, \bar{b}' \in A^k $, for some 
$h,k \geq 0$, $t$ is an $h+k$-ary term operation of $\mathbf A $,
and $\bar{a} \RR \bar{a}'$, $\bar{b} \SSS \bar{b}'$.
Here, $\bar{a} \RR \bar{a}'$ means $a_1 \RR a'_1$, $a_2 \RR a'_2$, \dots.

For $R, S, T, U$
compatible and reflexive relations,
let
\[
K( R, S; U; T) =
\left\{ 
(z,w) | 
\begin{vmatrix}   
x & y \cr z & w 
\end{vmatrix}  
\in M( R , S ), \ x \UU z, \ x \TT y,
\right\}. 
\] 

$K( R, S; 1; T) $ is denoted by
$ K( R, S; T)$ in Part I, and various commutators
for relations are constructed there from
$ K( R, S; T)$.

\begin{proposition}\label{hom}
For every variety $\v$,
$ K( R, S; U; T)$ and $ K( R, S; T)$
are monotone
 global operators on $\v$ satisfying
the homomorphism property.
\end{proposition}

\begin{proof}
We shall give the proof for $ K( R, S; U; T)$.
The proof for $ K( R, S; T)$ is similar and simpler 
(and, anyway, is the particular case $U=1$ constantly).
The operators are trivially monotone.

We now prove that $ K( R, S; U; T)$ satisfies the homomorphism property.
Suppose that ${\alg A}, {\alg  B} \in \v $,
$\phi: {\alg B} \to {\alg A} $,
$R, S, T, U \in \admb$,
 $x', w' \in {\alg A} $,
and 
$(x', w') \in \phi (K( R, S; U; T))$. 
We have to show that
$(x', w') \in K( \phi (R), \phi ( S); 
\brfr
 \phi (U); \phi  (T))$. 

That
$(x', w') \in \phi (K( R, S; U; T))$ means that 
$ {\alg B}$ has couples $(x_1,w_1), \dots, 
\brfr
(x_m, w_m) \in K( R, S; U; T)$
such that, in $ {\alg A} $,  $(x', w')$ belongs to the admissible and reflexive
relation generated by    
$( \phi (x_1), \phi (w_1)), \dots, ( \phi (x_m), \phi (w_m))$.
This means that $ {\alg A} $
has an $m$-ary  polynomial $p$ such that   
$x'=p( \phi (x_1), \dots,  \phi (x_m))$, and
$w'=p( \phi (w_1), \dots,  \phi (w_m))$.

Since $(x_i,w_i) \in K( R, S; U; T)$ for every $i=1,\dots m$, then,
by the definition of $ K( R, S; U; T)$, there are matrices
\[
\begin{vmatrix}
t_i( \bar{a}_i, \bar{b}_i)     & t_i( \bar{a}_i, \bar{b}_i') \cr
t_i( \bar{a}_i', \bar{b}_i)     & t_i( \bar{a}_i', \bar{b}_i') 
\end{vmatrix}   
\] 
such that 
$\bar{a}_i \RR \bar{a}_i'$, $\bar{b}_i \SSS \bar{b}_i'$,
$ t_i( \bar{a}_i, \bar{b}_i) \TT t_i( \bar{a}_i, \bar{b}_i') $,
$ t_i( \bar{a}_i, \bar{b}_i) \UU t_i( \bar{a}_i', \bar{b}_i)$, and 
$x_i=t_i( \bar{a}_i', \bar{b}_i) $, $ w_i=t_i( \bar{a}_i', \bar{b}_i')$.

For sake of brevity, let us write
$ t_i(\phi \bar{a}_i, \phi \bar{b}_i)$ in place of 
$ t_i(\phi (a_{i1}), \phi (a_{i2}), \dots, \brfr \phi (b_{i1}), \phi (b_{i2}), \dots)$.
The matrix
\[
\begin{vmatrix}
p\big(t_1( \phi \bar{a}_1, \phi \bar{b}_1), \dots, t_m( \phi \bar{a}_m, \phi \bar{b}_m)   \big)
 & p\big (t_1( \phi \bar{a}_1, \phi \bar{b}_1'), \dots, t_m( \phi \bar{a}_m, \phi \bar{b}_m') \big) \cr
p\big (t_1( \phi \bar{a}_1', \phi \bar{b}_1), \dots, t_m( \phi \bar{a}_m', \phi \bar{b}_m) \big)
  & p\big (t_1( \phi \bar{a}_1', \phi \bar{b}_1'), \dots, t_m( \phi \bar{a}_m', \phi \bar{b}_m') \big)
\end{vmatrix}
\] 
belongs to $ M( \phi (R) , \phi ( S) )$; moreover, since 
$ \phi $ is a homomorphism, the above matrix equals
\[
\begin{vmatrix}
p\big( \phi (t_1( \bar{a}_1, \bar{b}_1)), \dots, \phi (t_m( \bar{a}_m, \bar{b}_m))   \big)
 & p\big ( \phi ( t_1( \bar{a}_1, \bar{b}_1')), \dots, \phi (t_m( \bar{a}_m, \bar{b}_m')) \big) \cr
p\big ( \phi (t_1( \bar{a}_1', \bar{b}_1)), \dots, \phi (t_m( \bar{a}_m', \bar{b}_m)) \big)
  & p\big ( \phi (t_1( \bar{a}_1', \bar{b}_1')), \dots, \phi (t_m( \bar{a}_m', \bar{b}_m')) \big)
\end{vmatrix}   
\] 
Since 
$ t_i( \bar{a}_i, \bar{b}_i) \TT t_i( \bar{a}_i, \bar{b}_i') $,
and since $T$ is admissible,
the elements in the upper row of the matrix are $ \phi (T)$-related;
similarly, the elements 
in the left-hand column are $ \phi (U)$-related. 
Since $x_i=t_i( \bar{a}_i', \bar{b}_i) $, and $ w_i=t_i( \bar{a}_i', \bar{b}_i')$, then the second row of the matrix consists of the
elements
$ p( \phi (x_1), 
 \dots,  \phi (x_m))=x'$, and
$ p( \phi (w_1), \dots,  \phi (w_m)) =w'$.
This means that the matrix witnesses that 
$(x', w') \in K( \phi (R), \phi ( S);
 \brfr
 \phi (U);  
 \phi  (T))$. 
\end{proof}

\begin{proposition}\label{hom1}
If $\v$ is a variety,
 $F_1, F_2, \dots$ are $n$-ary
global operators on $\v$ for admissible and reflexive relations, 
and $F_1, F_2, \dots$ satisfy the homomorphism property, then so do
the following operators:

(i) $G_1(\bar{R})=F_1 (\bar{R}) \circ F_2 (\bar{R})$;

(ii) $G_2(\bar{R})=F_1 (\bar{R}) \cap F_2 (\bar{R})$;

(iii) $G_3(\bar{R})= \overline{F_1 (\bar{R}) \cup F_2 (\bar{R})}$;

(iv) $G_4(\bar{R})=F_1 (\bar{R}) + F_2 (\bar{R})$;

(v) $G_5(\bar{R})=\big(F_1 (\bar{R})\big)^*$;

(vi) $G_6(\bar{R})=
K\big(F_1 (\bar{R}),F_2 (\bar{R}); F_3 (\bar{R}); F_4 (\bar{R})\big)$.

(vii) $G_7(\bar{R})=
K\big(F_1 (\bar{R}),F_2 (\bar{R}); F_4 (\bar{R})\big)$.

If, in addition, $F_1, F_2, \dots$ are monotone, then so are
$G_1, \dots, G_7$.
\end{proposition}

\begin{proof}
The proof is similar to the proof of Proposition 
\ref{hom}, in many cases simpler. 
\end{proof}

Since, trivially,  the constant operators $F(\bar{R})=0$ and $ F(\bar{R})=1$, the operator $F(R)=R^-$, 
as well as 
the projections,
are monotone and
satisfy the homomorphism property, 
we get from Proposition \ref{hom1}(vii)(v)
that, for example, 
$F(R,S)=K(R,S;0)^*$
is monotone and satisfies the homomorphism
property.
Similarly, 
Proposition \ref{hom1}
can be iterated in order to show that 
all the commutators introduced in Part I
are monotone and satisfy the homomorphism property. 
Notice also that condition (vi) in Proposition \ref{hom1}
is more general than Proposition \ref{hom}: just take $n=4$, and let  each 
$F_i$ be the projection onto the $i$-th coordinate.

\end{document}